\documentclass{amsart}

\usepackage[utf8]{inputenc}
\usepackage[T1]{fontenc}
\usepackage{geometry}
\usepackage{amsmath,amssymb,amsthm}
\usepackage{graphicx}
\usepackage{hyperref}
\hypersetup{
    colorlinks=true,
    linkcolor=blue,
    citecolor=red,
    urlcolor=blue
}
\usepackage{cleveref} % For smart cross-referencing
\usepackage{booktabs}
\usepackage{siunitx} % For proper number and unit formatting

\newcommand{\Z}{\mathbb{Z}}

\newcommand{\alp}{\alpha}
\newcommand{\alpc}{\alpha^*}
\newcommand{\matM}{\mathbf{M}}
\newcommand{\vectS}{\mathbf{S}}
\newcommand{\states}{\mathcal{S}}
\newcommand{\forbidden}{\mathcal{F}}
\DeclareMathOperator{\specrad}{\rho}

\newtheorem{theorem}{Theorem}[section]
\newtheorem{definition}[theorem]{Definition}

\theoremstyle{remark}

\title[An Improved Lower Bound for the Stavskaya's Process]{An Improved Lower Bound for the Critical Parameter of the Stavskaya's Process via a Generalized Recurrent Method}
\author{Olivier Couronné}
\subjclass[2020]{60K35}

\address{Universit\'e Paris Nanterre, 200 avenue de la R\'epublique 92000 Nanterre, France, Modal'X, FP2M, CNRS FR 2036.}
\email{olivier.couronne@parisnanterre.fr}
\thanks{The author is supported by the Labex MME-DII funded by ANR,
reference ANR-11-LBX-0023-01, and this research has been conducted within the FP2M federation (CNRS FR 2036)}
\keywords{Particle Random Process;  Stavskaya’s Process}

\date{\today}

\begin{document}

\maketitle

\begin{abstract}
Stavskaya's process, a discrete-time version of the contact process on $\Z$, is known to exhibit a phase transition at a critical parameter $\alpc$ whose exact value remains an open problem. Recent work by Ramos et al. established a lower bound by linking the process's survival to the non-percolation of a dual contour. The probability of this contour was estimated using a recurrent method on a state space of weighted random walks with short-term memory. In this paper, we generalize and extend this method by systematically increasing the walk's memory and enriching the set of forbidden path sequences. By increasing the memory up to a length of 20 steps (corresponding to our  parameter $n=7$), we formulate the problem with a one-step transition matrix and numerically optimize its spectral radius. We thus establish the new lower bound $\alpc > 0.1370721$.
\end{abstract}

\section{Introduction}

The Stavskaya's process is a one-dimensional probabilistic cellular automaton, often considered the discrete-time analogue of the well-known contact process. The system evolves on the integer lattice $\Z$, with each site being in state 0 (vacant) or 1 (occupied). The dynamics are governed by a parameter $\alp \in [0,1]$, which represents the probability for an occupied site to become vacant.

The process exhibits a well-established phase transition at a critical parameter $\alpc\in(0,1)$. For $\alp > \alpc$, the process is ergodic and converges to the "all-0" configuration. For $\alp < \alpc$, a non-trivial stationary measure exists, and the process started from the "all-1" configuration survives with positive probability. While the existence of $\alpc$ is known, its exact value remains an open problem. Rigorous bounds were first established by Toom \cite{Toom1990}.

Ramos et al. \cite{Ramos2020} recently revisited this problem using a powerful duality argument. They mapped the survival of the process to the absence of an open contour in a dual graph. The probability of such a contour was bounded by summing the weights of all "nice paths" on this dual lattice. In their work, paths were deemed "nice" if they avoided a small set of short forbidden sequences. This approach, based on a recurrence relation for paths with a two-step memory, yielded a lower bound of $\alpc > 0.11$.

In this work, we build upon and generalize this recurrent method. The precision of the bound is fundamentally limited by the "memory" of the random walk used to model the contours. By increasing this memory and introducing more sophisticated forbidden patterns, we can capture the process's correlations more accurately, leading to a tighter bound. Our main contribution is to generalize the method as a function of an integer $n \ge 1$, where the walk's memory has length $3n-1$ and where we systematically forbid all primitive loops of length up to $3n$. This generalization unifies the approach and makes it adaptable to the available computational power.

Our main result, obtained through numerical optimization for the case $n=7$, is the following theorem.

\begin{theorem}
The critical parameter $\alpc$ of the one-dimensional Stavskaya's process satisfies
$$ \alpc > 0.1370721. $$
\end{theorem}

The remainder of this paper details the methodology. In \Cref{sec:method}, we define the generalized recurrent method, including the state space and the construction of the transition matrix as a function of the  parameter $n$. \Cref{sec:optimization} describes the numerical procedure used to find the optimal bound and presents the results for different values of $n$. We conclude in \Cref{sec:conclusion}.

\section{The Generalized Recurrent Method}
\label{sec:method}

We follow the framework of \cite{Ramos2020}, which itself builds on the method of Toom \cite{Toom1990}, where the probability of the process dying out is bounded by the probability of forming a specific type of contour on a dual lattice. 
This probability is bounded by the sum of weights of all valid contour paths.

In the original method, this sum of path weights is analyzed using a system of generating functions. 
They defined quantities $S_r(k)$ representing the total weight of all valid paths of length $k$ that terminate with a specific move $r$. 
The weights are functions of the process parameter $\alp$ and two auxiliary real parameters, $p, q \ge 1$. 
The analysis then relies on deriving a system of linear recurrence relations that connect these sums at different lengths (e.g., expressing sums at length $k+2$ in terms of sums at length $k$).
The convergence of the total sum is then determined by the properties of this system of recurrences.

Our approach generalizes this concept. Instead of tracking the sums based on the last step, we define a comprehensive state space where each state represents a much longer path history. 
This allows us to reframe the problem: the complex, multi-step recurrence relations of the original method become a simple one-step linear transformation governed by a large transition matrix $\matM_n$. 
The condition for convergence then translates to the spectral radius of this matrix being less than one.
The core of our new approach is to systematically generalize and refine the definition of a "valid" path history to construct this matrix.

\subsection{Paths, Weights, and Primitive Loops}
On the dual graph, a path is a sequence of steps. The three possible step types have a weight that depends on the parameters $p, q \ge 1$ and $\alp$:
\begin{itemize}
    \item Step 1 (displacement $(-1,-1)$): weight $p^{-1}q^{-1}$
    \item Step 2 (displacement $(+2,0)$): weight $\alp p^{2}$
    \item Step 3 (displacement $(-1,+1)$): weight $qp^{-1}$
\end{itemize}

A key insight of the contour method is that not all sequences of steps form valid contours. 
Following Toom \cite{Toom1990} and Ramos et al. \cite{Ramos2020}, we identify two types of forbidden patterns.

First, the sequences $(1,3)$ and $(3,1)$ are forbidden. 

Second, we exclude all closed loops. Given the vector displacements of the three step types, a path returns to its origin—thereby forming a closed loop—if and only if it is "balanced," meaning it contains an equal number of steps of type 1, 2, and 3.

\begin{definition}
A \textbf{loop of order $k$} is a path of length $3k$ containing exactly $k$ instances of each step type $\{1, 2, 3\}$. A loop of order $k$ is said to be \textbf{primitive} if it does not contain a primitive loop of order $j < k$ as a contiguous subpath.

The set of forbidden patterns $\forbidden_n$ for a  level $n$ consists of:
\begin{enumerate}
    \item The patterns $(1,3)$ and $(3,1)$.
    \item All primitive balanced loops of order $k$, for $1 \le k \le n$.
\end{enumerate}
\end{definition}

For example, for $n=1$, $\forbidden_1$ contains the degenerate loops $(1,3)$ and $(3,1)$, as well as the primitive loops of order 1: $(1,2,3)$ and $(3,2,1)$. This set of patterns corresponds exactly to the forbidden sequences used in the analysis by Ramos et al. \cite{Ramos2020}.

For $n=2$, we add the 2 primitive loops of order 2, that is $(1,1,2,2,3,3)$ and $(3,3,2,2,1,1)$.

\subsection{State Space and Transition Matrix}

For a  level $n$, we define a state space based on a memory of the last $L = 3n-1$ steps.

\begin{definition}[State Space]
For a  level $n$, the state space $\states_n$ is the set of all paths of length $L=3n-1$ that do not contain any primitive loop from $\forbidden_{n-1}$ as a subpath.
\end{definition}

Let $\vectS(k)$ be a column vector where each entry $S_w(k)$ is the sum of weights of all valid paths of length $k$ ending in the state $w \in \states_n$. A path of length $k+1$ is formed by appending a new move $j \in \{1, 2, 3\}$ to a valid path of length $k$. This transition is allowed only if the resulting sequence at the tail of the path does not create any forbidden pattern.

The evolution of the weighted path counts is described by a one-step linear transformation:
\begin{equation}
\vectS(k+1) = \matM_n(p, q, \alp) \vectS(k),
\end{equation}
where $\matM_n(p, q, \alp)$ is the transition matrix of size $|\states_n| \times |\states_n|$. 
To make the dependence on the parameters explicit, we can decompose $\matM_n$ into three base matrices. 
Each base matrix, $\mathbf{M}^{(j)}$, exclusively represents the valid transitions ending with the specific move $j \in \{1, 2, 3\}$:
\begin{equation}
\matM_n(p, q, \alp) = \frac{1}{pq}\mathbf{M}^{(1)} + (\alp p^2)\mathbf{M}^{(2)} + \frac{q}{p}\mathbf{M}^{(3)}.
\label{eq:matrix_decomposition}
\end{equation}
Here, $\mathbf{M}^{(j)}$ for $j \in \{1, 2, 3\}$ are sparse matrices with entries of 0 or 1.

An entry $M_{w', w}$ is non-zero if and only if state $w'$ can be obtained from state $w$ by appending a valid move $j$. 
More formally, if $w=(w_1, \dots, w_L)$, a transition to $w'=(w_2, \dots, w_L, j)$ is valid if the extended path $(w_1, \dots, w_L, j)$ of length $3n$ does not contain any pattern from $\forbidden_n$ as a suffix. 
If the transition is valid, $M_{w', w}$ takes the value of the weight of the new move $j$.

\section{Numerical Optimization and Results}
\label{sec:optimization}

The total probability of a contour existing converges if the spectral radius of the transition matrix, $\specrad(\matM_n)$, is less than 1. Our goal is to find the largest value of $\alp$ for which we can find parameters $p, q \ge 1$ that satisfy this condition.
\begin{equation}
\alpc \ge \sup_{n \ge 1} \sup_{\substack{p \ge 1 \\ q \ge 1}} \sup\left\{\alp \in [0,1] \mid \specrad\bigl(\matM_n(p,q,\alp)\bigr) < 1 \right\}.
\end{equation}

Our numerical explorations indicated that the maximum is achieved along the line $q=1$. 
The optimization was therefore reduced to a one-dimensional search over $p$.
 
For a fixed $p$ (and $q=1$), our goal is to find the supremum of the set of all $\alpha$ values for which the condition $\specrad(\matM_n(p, 1, \alpha)) < 1$ is satisfied. We denote this supremum as $\alpha_c(p)$.

Numerically, we approximate $\alpha_c(p)$ using a bisection search. The search iteratively narrows down an interval $[\text{low}, \text{high}]$ for $\alpha$. At each step, a midpoint $\text{mid}$ is tested. If $\specrad(\matM_n(p, 1, \text{mid})) < 1$, we know that $\text{mid}$ is a valid lower bound, so we set $\text{low} = \text{mid}$. Otherwise, we set $\text{high} = \text{mid}$. This process converges to a highly accurate estimate of $\alpha_c(p)$ while ensuring that the final $\text{low}$ value always corresponds to a spectral radius less than 1.

The best lower bound for $\alpc$ is then the maximum value found over our search space for $p$:
$$ \alpc > \sup_{p \ge 1} \alpha_c(p). $$

The numerical procedure for a fixed $n$ is as follows:
\begin{enumerate}
    \item Generate the set of primitive loops $\forbidden_n$.
    \item Construct the state space $\states_n$ and the base transition matrices.
   \item For each value of $p$ in a search grid (with $q=1$):
    \begin{enumerate}
        \item Approximate the supremum $\alpha_c(p)$ of the set $\{\alpha \mid \specrad(\matM_n(p, 1, \alpha)) < 1\}$ using a bisection search with a tolerance of $10^{-10}$.
    \end{enumerate}
    \item The best lower bound for $\alpc$ at level $n$ is the maximum value found over the grid, $\max_p \alpha_c(p)$.
\end{enumerate}

The results of this procedure for increasing values of $n$ are presented in \Cref{tab:results}. 
We observe of course an improvement in the lower bound, at the cost of an exponential growth in the size of the state space.

\begin{table}[htbp]
    \centering
    \caption{Lower bounds for $\alpc$ as a function of the  level $n$.}
    \label{tab:results}
    \begin{tabular}{crrS[table-format=1.4]S[table-format=1.8]}
        \toprule
        \textbf{n} & \textbf{Primitive Loops} & \textbf{State Space Size} & {\textbf{$p_{opt}$}} & {\textbf{Lower Bound for $\alpc$}} \\
        \midrule
1 & 4 & 7 &1.464  & 0.125 \\
2 & 6 & 73 & 1.44 &  0.13101966\\
3 & 12 & 759 & 1.43 & 0.13358660 \\
4 & 36 &7,859  &1.424  &  0.13502855\\
5 & 146 &81,231  &1.42  & 0.13595342 \\
6 & 694 & 839,009 &1.417  &0.13659747  \\
7 & 3584 &8,663,071  &1.415  &0.13707211  \\
        \bottomrule
    \end{tabular}
    \flushleft
    \footnotesize{\textbf{Note:} The "State Space Size" column corresponds to $|\states_n|$. $p_{opt}$ is the value of $p$ that maximizes the bound for $q=1$.}
\end{table}

The search for $n=7$ yields our main result. The optimal value was found for parameters $p = 1.415$ and $q=1.0$.

\section{Conclusion}
\label{sec:conclusion}

By generalizing the recurrent method of Ramos et al., we have developed a systematic framework for improving the lower bound of the critical parameter of the Stavskaya's process. By increasing the memory of the random walk and refining the set of forbidden paths through the algorithmic exclusion of primitive loops, we have established a new lower bound of $\alpc > 0.1370721$.

This work demonstrates the power of the dual contour method when combined with a more detailed state representation and numerical optimization. The framework is general and can be applied with even larger memory lengths ($n > 7$), which would yield further improvements. However, the exponential growth in the size of the state space, visible in \Cref{tab:results}, constitutes a major computational bottleneck, suggesting that future gains from this approach alone will be increasingly costly to obtain.

\end{document}